\documentclass{article}
\usepackage{amsmath,amsfonts,latexsym}
  \usepackage{paralist}
  \usepackage{graphics} %% add this and next lines if pictures should be in esp format
  \usepackage{epsfig} %For pictures: screened artwork should be set up with an 85 or 100 line screen
\usepackage{graphicx}  \usepackage{epstopdf}%This is to transfer .eps figure to .pdf figure; please compile your paper using PDFLeTex or PDFTeXify.
 \usepackage[colorlinks=true]{hyperref}
\hypersetup{urlcolor=blue, citecolor=red}

\newtheorem{Theorem}{Theorem}[section]

\newtheorem{Proposition}{Proposition}[section]

\newtheorem{Corollary}{Corollary}[section]
\newenvironment{Proofc}[1]{\smallskip\par\noindent\textsc{#1}\quad}%
  {\hfill$\Box$\bigskip\par}
\newenvironment{Proof}{\begin{Proofc}{Proof}}{\end{Proofc}}

\newtheorem{Remark}{Remark}[section]

\def\d{\delta}

\def\g{\gamma}
\def\G{\Gamma}
\def\l{\lambda}

\def\r{\rho}

\def\z{\zeta}
\def\pd{\partial}
\def\half{\frac{1}{2}}
\newcommand{\cE}{{\mathcal E}}

\newcommand{\cG}{{\mathcal G}}

\newcommand{\cN}{{\mathcal N}}

\newcommand{\cV}{{\mathcal V}}
\newcommand{\R}{{\mathbb R}}
\newcommand{\N}{{\mathbb N}}

\newcommand{\ud}{{\mathrm{d}}}

\def\pd{\partial}

\begin{document}
\title{A discrete Hughes' model for pedestrian flow on graphs}
\author{Fabio Camilli\footnotemark[1], \,Adriano Festa\footnotemark[2], \, Silvia Tozza\footnotemark[3]}
\date{version: \today}
\maketitle

\footnotetext[1]
{
 Dip. di Scienze di Base e Applicate per l'Ingegneria,  ``Sapienza" Universit{\`a}  di Roma,
 via Scarpa 16, 00161 Roma, Italy, ({\tt e-mail:camilli@sbai.uniroma1.it})
}

\footnotetext[2]
{RICAM,  Austrian Academy of Sciences (\"OAW)
Altenbergerstr. 69, 4040 Linz, Austria. ({\tt e-mail: a.festa@oeaw.ac.at})
}

\footnotetext[3]
{
Dipartimento di Matematica,  ``Sapienza" Universit{\`a}  di Roma,
P.le Aldo Moro, 5 - 00185 Roma, Italy
({\tt e-mail: silvia.tozza@uniroma1.it})
}
\begin{abstract}
In this paper, we introduce a discrete time-finite state  model for pedestrian flow on a graph
in the spirit  of the  Hughes dynamic continuum model. The pedestrians, represented by a  density function,
move  on the graph choosing a route to minimize the instantaneous travel cost to the destination.
The  density is  governed by a  conservation law while the minimization principle
is described by a graph eikonal equation. We show that the model is well posed and we implement
some numerical examples  to demonstrate the validity  of the proposed  model.
\end{abstract}
 \begin{description}
%\item [\textbf{ AMS subject classification}:  ]
%%
  \item [\textbf{Keywords}: ] Pedestrian flow, Hughes' model, conservation law, graph  eikonal equation.
 \end{description}
%
%
%
%%%%%%%%%%%%%%%
%%%        Introduzione        %
%%%%%%%%%%%%%%%
%
\section{Introduction}\label{intro}
In \cite{h00}, Hughes introduced a by now classical  fluidodynamical  model   to study the motion  of a large human crowd
 (see also \cite{hw}, \cite{h02}, \cite{h03},  \cite{tcp} and \cite{bd} for a review).

The crowd is treated as a ``thinking fluid'' and it moves at the maximum speed towards   a common destination or goal, taking also into account the environmental conditions. In fact, people prefer  to avoid crowded regions and this assumption is incorporated in a  potential field  which gives the direction of the motion. The potential field is described by  the solution of an eikonal equation giving the optimal  paths to the destination, integrated with respect to a cost proportional to the local crowd density. Hence, for each instant of time, an individual looks at the global configuration of the crowd and updates his/her direction to the exit trying to avoid the crowded, motion expensive regions.

Many situations related to pedestrian motion,  for example the study of a crowd escaping  from a building,  can be
described as the problem of finding the shortest path in a network.  A model to simulate the behavior of pedestrian motion
on networks is therefore important for designers to analyze the performance results in terms of number of nodes and  connectivity of the environment.

In this paper we  introduce  both a  model for pedestrian motion on networks which on one side can be viewed as a discrete time-finite state analogous of the Hughes'  model, and a numerical discretization of the Hughes' system  defined on a network. The system is composed by a graph eikonal equation  where the cost depends on the density of the population, and a graph conservation law  which governs the  evolution   of the population. We show that, under some natural assumptions on the flux, the  graph Hughes' system is well-posed for any time $n\in \N$. Moreover it shares with the corresponding continuous model some qualitative properties, like e.g. the interpretation of the solution of the graph eikonal equation as a distance from the boundary.

The model here described bears some resemblance to the \emph{rational behavior model} studied  in \cite{cp}. People know
at each  time the global distribution of the population on the network and therefore, they accordingly modify their strategy to reach the exit.
In \cite{cp} this behavior is obtained by   introducing an optimization problem at the junctions, whereas here the optimal strategy is given by
 the solution of the eikonal equation.

%We first analyze a discrete eikonal equation on the graph where the cost depends on the distribution of the population and which becomes unbounded %\textcolor{blue}{ for the distribution going to 1}; we prove existence and uniqueness of the solution, giving also a representation formula for the solution in %terms of a distance from the boundary. Then, we analyze the  conservation law  governing the distribution of the population and we show that, under some %natural assumptions on the flux, \textcolor{blue}{ if the initial density $\bar\r$ is strictly less than 1,  this property is preserved for any time}. It %follows that the  discrete eikonal equation  is well-posed and, as a consequence, there is a solution of the graph Hughes' model for any time $n\in \N$.
%Hence, in contrast with the continuous Hughes' model, its discrete analogous is well-posed, i.e. its solution is well defined  for all times $n\in\N$. \\

We also present an algorithm for the solution of the discrete problem. %described.
%The final section is devoted to an algorithm for the solution of the discrete problem.
We consider several examples  and we show that the discrete Hughes' system  captures the natural behavior of the crowd.\par
The paper is organized as follows. In Section \ref{Sec1}, we recall the Hughes' model and we formulate its analogous on a network. In Section \ref{Sec2},
we derive the discrete Hughes' system on a graph and in Section \ref{Sec3} we prove the well-posedness of this problem. Section \ref{Sec4} is devoted
to the algorithm for the solution of the problem and to numerical experiments aimed at confirming the validity of the proposed model.
We conclude the paper  with final comments and remarks reported in Section \ref{Sec5}.

%%%%%%%%%%%%
%          %
%%%%%%%%%%%%
\section{The Hughes' model}\label{Sec1}
From a mathematical point of view the Hughes' model consists in the following system
\begin{align}
&\partial_t \rho-\hbox{div}(\rho f^2(\rho)\nabla u)=0,&x\in\Omega, \label{hughes1}  \\[4pt]
&|\nabla u|=\frac{1}{f(\rho)}, &x\in\Omega, \label{hughes2}
\end{align}
where $\Omega$ is a bounded domain of $\R^n$ and    $\rho$, which takes value in $ [0,1]$,  is a density field representing  the concentration  of the pedestrian at $(x,t)$.
The problem is completed with  some boundary conditions: the initial configuration of the mass
\begin{equation*}\label{initial_mass}
\rho(x,0)=\bar \rho(x), \quad x\in \Omega,
\end{equation*}
and a no-flux condition on the boundary
\begin{equation}\label{no_flux_BC}
\rho f^2(\rho)\nabla u=0, \quad x\in\partial \Omega,
\end{equation}
for the continuity equation \eqref{hughes1}; the Dirichlet boundary condition
\begin{equation*}\label{Dirichlet_BC}
u(x)=0, \quad x\in\partial \Omega,
\end{equation*}
for the eikonal equation  \eqref{hughes2}.
The  function $f(\r)$ is typically given by $1-\r$, hence    the  (absolute value of the)  flux   is given, as   in many other models related
to  pedestrian and vehicular motion, by $g(\r)=\r v_{max} (1-\r)$. The term
$v_{max}$, which can be always  normalized to 1, is  an maximal speed at which an agent would travel in ideal environment conditions and the term $\rho(1-\rho)$, the  so called fundamental diagram (see \cite{gp}), indicates that the velocity of a   pedestrian is proportional
to the crowd density (recall that $\r\in [0,1]$).\\
The direction of the motion is given by the potential term  $\nabla u$. If  the cost $1/f(\rho(t))$ is bounded, i.e. if  $\rho(t)<1$,
 the solution of the  eikonal equation \eqref{hughes2} at time $t$ is given by
\begin{equation}\label{dist_boundary}
u(x)=\inf\{d_t(x,y):\, y\in \partial \Omega\}
\end{equation}
where  the   distance function $d_t:\Omega\times\Omega\to\R^+$   is defined as follows
\[
d_t(x,y)=\inf\left\{\int_0^S\frac{1}{f(\rho(\xi(s),t))} \ud s:\,S>0,\,\xi\in G^S_{x,y}\right\},
\]
with $G^S_{x,y}$  the set of the absolute continuous curves in $\Omega$ such that $\xi(0)=y$, $\xi(S)=x$ and $|\dot\xi(s)|=1$ a.e. in $[0,S]$.
The term $1/f(\rho(t))$  can be interpreted as the current  cost associated to the curve $\xi(s)$ and the solution $u$ of \eqref{hughes2} selects the curves which minimize  the cost for reaching the boundary. Hence, people  are directed towards the boundary trying to avoid crowded regions. \\
Existence of a solution  to \eqref{hughes1}-\eqref{hughes2} is still an  open problem, the main difficulty is
given by the possible concentration of population for some $t$ which results in the blow-up of the cost term $1/f(\rho(t))$. Partial results are available only  in the one-dimensional case (see \cite{amad12,amad14,dmpw});   note that in $\R$, since  $|\pd_x u| = \pd_x u \,\hbox{sign}(\pd_x u)$, the first equation in \eqref{hughes1}   simplifies as
\begin{equation*}\label{hughes1d}
\partial_t \rho-\pd_x\left(\rho f(\rho) \,\hbox{sign}(\pd_x u )\right)=0
\end{equation*}
and   the solution of the eikonal equation admits an almost explicit representation formula.
Moreover, there have been several approaches which regularize the flux function in order to obtain a well-posed problem. In \cite{dmpw}, a regularization of the eikonal equation \eqref{hughes2} has been introduced in order to avoid the possible  blow-up of $|\nabla u|$, leading to the system
\begin{equation}\label{hughesreg}
\begin{split}
&\partial_t \rho-\hbox{div}(\rho f^2(\rho)\nabla u)=0,  \\[4pt]
&-\epsilon  \Delta u  +| \nabla u|^2=\frac{1}{(f(\rho)+\delta)^2},
\end{split}
\end{equation}
for some $\epsilon$, $\delta>0$.
 Also for this system, existence and uniqueness of a weak solution have been obtained in the one dimensional case.
%%%%%%%%%%%%
\subsection{The Hughes' model on networks}
In the recent years,  the theory of weak solutions for  conservation laws and of viscosity solutions  for  Hamilton-Jacobi equations % (see \cite{cm})
have been extended to the case of networks (see \cite{gp} and \cite{cm}, respectively), imposing appropriate transition conditions at the vertices. \par
A network  $\mathcal N=(\cV,\cE)$     is a finite  collection of points $\cV:=\{v_i\}_{i\in \mathcal I}$ in $\R^n$ connected by continuous, non self-intersecting edges $\cE:=\{e_j\}_{j\in \mathcal J}$.
We define $N:=|\cV|, M := |\cE|$ and $\mathcal I:=\{1,\dots,N\}$, $\mathcal J:=\{1,\dots,M\}$.
Each edge $e_j\in \cE$ is  parametrized by a smooth function $\pi_j:[0,l_j]\to\R^n,\, l_j>0$. Given $v_i\in \cV$,
$Inc_i:=\{j\in \mathcal J:v_i\in e_j \}$ denotes the set of edges branching out from $v_i$, and we denote by $d_{v_i}:=|Inc_i|$ the degree of $v_i$. A vertex $v_i$ is called a  boundary vertex if $d_{v_i}=1$ and  $\pd \cN$ denotes the set of boundary vertices. \\
For a function $u:  \cN\to\R$,    $u_j:[0,l_j]\to \R$ is the restriction of $u$ to $e_j$, i.e. $u(x)=u_j(y)$ for $x\in e_j$, $y=\pi_j^{-1}(x)$, and     $\pd_j u(v_i)$  is the  oriented derivative of $u$ at $v_i$ along the arc $e_j$ defined   by
%\[\pd_j u (v_i)=\lim_{x\in e_j,\,x\to v_i}\frac{u(x)-u(v_i)}{|\pi^{-1}_j(x)-\pi^{-1}_j(v_i)|}\]
\begin{equation*} %\[
\pd_j u (v_i)=\left\{
               \begin{array}{ll}
                 \lim\limits_{h\to 0^+}(u_j(h)-u_j(0))/h, & \hbox{if $v_i=\pi_j(0)$,} \\
                 \lim\limits_{h\to 0^+}(u_j(l_j-h)-u_j(l_j))/h, & \hbox{if $v_i=\pi_j(l_j)$.}
               \end{array}
             \right.
\end{equation*} %\]
The Hughes' system on the network $\cN$ is given by
\begin{eqnarray}\label{hughesnetwork}
 \left\{
 \begin{array}{ll}
  \partial_t\rho_j(x,t)-\pd_x\left(\rho_j(x,t) f(\rho_j(x,t))\,\hbox{sign}(\pd_x u_j )\right)  =0\quad& x\in e_j \,, t>0 \,, j\in \mathcal J,\\[4pt]
| \pd_x u_j | =\frac{1}{f(\rho_j(x,t))}\ &x\in e_j \, t>0 \,, \ j\in \mathcal J,\\ [4pt]
 \sum\limits_{j\in  Inc_i}\rho_j(v_i,t) f(\rho_j(v_i,t))\,\hbox{sign}(\pd_j u(v_i) )=0\,,& \ t>0\,,\ i\in \mathcal I,\\[4pt]
  u_j(v_i)=u_k(v_i)& j,k\in Inc_i\,,\ i\in \mathcal I,\\[4pt]
\rho_j (x,0)=\bar \rho_j (x)\,,&  x\in\cN, \,  j\in \mathcal J,\\[4pt]
u_j(x)=0\,,&  x\in\pd\cN\,,\ j\in \mathcal J,
 \end{array}
 \right.
\end{eqnarray}
where the derivatives $\pd_x u_j$ must be interpreted as derivatives with respect to $y=\pi_j^{-1}(x)$, which parametrizes the edge $e_j$. \\
The system \eqref{hughesnetwork}   is formally equivalent to   $M$ scalar Hughes' systems defined on the edges    coupled via   the transition conditions at the internal vertices where we require  for the density $\rho$ the  conservation of the flux, while   for  $u$ the continuity. We recall that if we consider the two equations separately  the previous transition conditions are sufficient
to get existence (and for the Hamilton-Jacobi equation also uniqueness) of an entropy  solution for the   conservation law and of the viscosity solution
for the  Hamilton-Jacobi equation. In principle these two approaches could be combined to study the Hughes' system on networks, but it seems very difficult   since  even in the Euclidean case the existence of a solution is still an open problem. Hence, in Section \ref{Sec2} we will consider a discretization of the system
\eqref{hughesnetwork} obtaining a well-posed discrete version of the Hughes' system. Finally, we observe that it is possible to consider   the following network version of the regularized  system \eqref{hughesreg}
\begin{eqnarray*}
 \left\{
 \begin{array}{ll}
  \partial_t\rho_j(x,t)-\pd_x\left(\rho_j(x,t) f(\rho_j(x,t))\,\hbox{sign}(\pd_x u_j )\right)  =0\quad& x\in e_j, \, t>0 \,,\ j\in \mathcal J,\\[4pt]
-\epsilon \pd_{xx} u_j+| \pd_x u_j |^2 =\frac{1}{(f(\rho_j(x,t))+\delta)^2}\ &x\in e_j\, , \, \ j\in \mathcal J,\\ [4pt]
 \sum\limits_{j\in  Inc_i}\rho_j(v_i,t) f(\rho_j(v_i,t))\,\hbox{sign}(\pd_j u(v_i) )=0\,,& \ t>0\,,\ i\in \mathcal I,\\[4pt]
  u_j(v_i)=u_k(v_i)& j,k\in Inc_i\,,\ i\in \mathcal I,\\[4pt]
  \sum\limits_{j\in  Inc_i}\epsilon\pd_j u(v_i) =0\,,&\ i\in \mathcal I,\\[4pt]
\rho_j (x,0)=\bar \rho_j (x)\,,&  x\in\cN\,, \, j\in \mathcal J,\\[4pt]
u_j(x)=0\,,&  x\in\pd\cN\,,\ j\in \mathcal J,
 \end{array}
 \right.
\end{eqnarray*}
where, as stated above, the derivatives $\pd_x u_j$ and $\pd_{xx} u_j$ must be interpreted as derivatives with respect to $y=\pi_j^{-1}(x)$, which parametrizes the edge $e_j$.
Note that, with respect to the first order system, a Kirchhoff law for $u$ at the internal vertices has to be added.
\section{The Hughes' model on graph}\label{Sec2}
In this section, after a preliminary paragraph that introduces our notation on graph, we focus on the discrete Hughes' system and its interpretation as a discrete-time finite state model for pedestrian flow on a graph.
\subsection{Basic notations}
Let us consider a weighted  graph $\G=(V,E,w)$ where $V$ denotes the set of vertices,  $E\subset V\times V$  the set of the edges  and $w:V\times V\to \R$ the weights of the edges  with $w(x,y)>0$  if $(x,y)\in E$ and $w(x,y) = 0$ otherwise. In what follows we will use the notation $w_{xy} := w(x,y)$ for the sake of simplicity.
 %$w_{xy}>0$  if $(x,y)\in E$ and $w_{xy}=0$ otherwise.
 The weight $w_{xy}$ is a parameter  which takes into account several properties of the edge $(x,y)$
 such as length, capacity, velocity (small weights correspond  to better connection   between $x$ and $y$).
The graph is assumed to  be finite, simple, connected and undirected, hence $w_{xy}=w_{yx}$ for any $(x,y)\in E$.\\
We   set   $y\sim x$ if $(x,y)\in E$ and we  denote by $I(x)=\{y\in V:\, y\sim x\}$ the set of the neighbors of $x$  and
 by  $|I(x)|$  the degree of $x$, i.e. the cardinality of the set $I(x)$. We set
\begin{equation}\label{maxdeg}
    D=\max_{x\in V} |I(x)|.
\end{equation}
We denote by $V_b\subset V$  the set of  the boundary vertices of the graph
 and by $V_0=V\setminus V_b$   the set of the internal vertices.\\
A path connecting  $x$ to $y$  is given  by a finite subset  $\g=\{x_0=x,x_1,\dots, x_N=y\}$  of $V$ such that $x_k\sim x_{k+1}$, $k=0,\dots,N-1$.
We denote by $\cG_{xy}$ the set of the paths $\g$ connecting $x$ to $y$.
The geodetic distance between two adjacent vertices $x,y$ is
\[d(x,y)=w_{xy}\]
 whereas for two arbitrary vertices  $x,y\in V$ we define
\begin{equation}\label{geo}
d (x,y):=\min\{d(x_0,x_1)+d(x_1,x_2)+\dots+d(x_{N-1},x_N) \},
\end{equation}
where the minimum  is taken over all the finite paths $\g\in   \cG_{xy}$.
%%%%%%%%%%%%%%%%%%%%%%%%%%%%%
\subsection{The discrete Hughes' model}
Given a weighted  graph $\G$, we consider the following discrete Hughes' system for a $ n\in\N$
\begin{equation}\label{hughesD}
    \left\{ \begin{array}{ll}
       \r^{n+1}(x)= \r^{n}(x)-\sum\limits_{y\sim x}\l h^n_{yx}\cdot\textrm{sgn}(u^n(y)-u^n(x)),\quad &  x,y\in V, \\[4pt]
       \max\limits_{y\sim x}\left\{-\frac{u^n(y)-u^n(x)}{w_{yx}}- \frac{1}{1-\r^n(y)}\right\}=0, &x\in V_0, \, y\in V, \\[4pt]
        \r^0(x)=\bar\r(x), & x\in V,\\
        u^n(x)=0, &x\in V_b,
     \end{array}
     \right.
\end{equation}
where $h^n_{yx}$ denotes the flux between $x$ and $y$ and $\l$ is a positive constant.
To define the flux $h^n_{yx}$ in \eqref{hughesD} we consider a function  $h$   satisfying
\begin{eqnarray}
   & h(0,0)=h(1,1)=0,		\label{flux0}\\
   & m_-(v)\le h_v(v,u)\le 0\le h_u(v,u)\le m_+(u),	 \label{flux1}
\end{eqnarray}
for a continuous function $m:\R \to\R$, where $a_-=\min(a,0)$, $a_+=\max (a,0)$. We set
\begin{equation}\label{flux}
   h^n_{yx}=\left\{
              \begin{array}{ll}
               h(\r^n(y),\r^n(x)), & \hbox{if $\d^n_{yx}=1$,} \\[4pt]
                h(\r^n(x),\r^n(y)), & \hbox{if $\d^n_{yx}=-1$,}
              \end{array}
            \right.
\end{equation}
where $\d^n_{yx} = \textrm{sgn}(u^n(y)-u^n(x))$ with $u^n$  given by the second equation in \eqref{hughesD}. \\
In order to give specific examples of $h$,  we consider flux functions     which are consistent with $g(\r)=\r(1-\r)$, i.e.   $h(\r,\r)=g(\r)$ for $\r\in\R$, such as the Lax-Friedrichs
flux
\begin{equation}\label{LF}
      h(\r^n(y),\r^n(x))=\half\Big(\r^n(y)(1-\r^n(y))+\r^n(x)(1-\r^n(x))\Big)-\frac{1}{\l}(\r^n(y)-\r^n(x)).
\end{equation}
Other examples of  flux  verifying \eqref{flux0}-\eqref{flux1} are given by the \emph{Godunov flux}
\begin{equation*}\label{God}
h(\r^n(y),\r^n(x))=\left\{
         \begin{array}{ll}
           \min\limits_{[\r^n(x),\r^n(y)]} g(\r), & \hbox{if $\r^n(x)\le \r^n(y)$,} \\[4pt]
           \max\limits_{[\r^n(y),\r^n(x)]} g(\r), & \hbox{if $\r^n(x)\ge \r^n(y)$,}
         \end{array}
       \right.
    \end{equation*}
and by the   \emph{Engquist-Osher flux}
\begin{equation}\label{EO}
      h(\r^n(y),\r^n(x))=\half\Big(\r^n(y)(1-\r^n(y))+\r^n(x)(1-\r^n(x))\Big)-\half\int_{\r^n(x)}^{\r^n(y)} |g'(\r)|\ud\r.
\end{equation}
In all the previous examples $h$ is consistent with $g$  and satisfies \eqref{flux0}-\eqref{flux1} with $m(v)=g'(v)$. Note that
 the discrete conservation law in \eqref{hughesD}    coincides with the numerical scheme for \eqref{hughes1} introduced in \cite{t} if   the graph $\G$
is given by the discretization  points of an interval $[a,b]$ and $h$ is consistent with $g(\r)=\r(1-\r)$
(see also \cite{hw}). In the case of a network,  it corresponds to discretize  the conservation law \eqref{hughes1} inside the edges and to impose the conservation of the flux at the vertices. We also refer to \cite{bc},\cite{gp} for different numerical discretizations of conservation law on networks in the framework of vehicular traffic motion.
Concerning the Eikonal equation \eqref{hughes2}, we recall that approximations of Hamilton-Jacobi equations on networks are discussed in  \cite{clm}  for finite differences and in \cite{cfs} for semi-Lagrangian schemes.\par
The system \eqref{hughesD} has been introduced as a discretization of the continuous problem \eqref{hughes1}-\eqref{hughes2},
but nevertheless it inherits some dynamical properties of the original model and it can be interpreted as a discrete-time finite state model for the flow of pedestrians on a graph in the following way.
At the initial time, there is a continuum of indistinguishable players distributed on the vertices of the graph $\G$
according to a density function  $\bar\r$, where $\bar\r(x)$ represents  the crowd   at the vertex $x$.
As in the original Hughes' model \cite{h00},  we assume that  people have a complete knowledge of the environment and they  choose the minimum distance
path to their destination $V_b$, but they have difficult to move  against the flow which is proportional to the local crowd density.
Each vertex $x\in V$ represents a point at which people can choose which route $(x,y)\in E$ to take and the subset $V_b\subset V$   represents the goal, e.g.   the exit of the environment.  \\
The term  $\textrm{sgn}(u^n(y)-u^n(x))$ (with $\textrm{sgn}(0)=0$),
where $u^n:V\to\R$ is the solution of the graph eikonal equation in \eqref{hughesD} at the time $n$,  gives the direction of the flow.
The function $u^n$,  see \eqref{dist_boundaryD} for more details, is given by the distance function from the boundary    integrated along  the cost
 $1/(1-\r^n)$ which penalizes  the vertices with high density $\r^n$.  \\
The pedestrians having reached  a vertex $x\in V$ at a given time $n$ are stirred
to the arc which minimizes the distance $u^n$ from the boundary.
Moreover, if $x\in V_b$ then $u^n(x)=0$  and  $\d^n_{yx}>0$ for each $y\in V_0$ such that $(x,y)\in E$. Hence, the flow is always in the direction  of the boundary, i.e.   the individuals having already reached the destination cannot reenter  in the environment.\par
If we consider the flow $h$ e.g. given by \eqref{LF}, the first of the two terms  says that the velocity of the pedestrians along the arc $(x,y)$ is given by the average of the velocity at $x$ and $y$; the second term   is  a small viscosity term  given by the graph laplacian
\begin{equation} \label{diff}
-\frac{1}{\l}\sum\limits_{y\sim x}(\r^n(y)-\r^n(x))\d^n_{yx}
\end{equation}
which can be interpreted as a  a stochastic perturbation of  the flux at $x$.
\section{Analytical results for the  Hughes'  model  on graph}\label{Sec3}
In this section we prove existence and uniqueness of the solution of the discrete Hughes' model. We  study separately the discrete eikonal equation
and the discrete conservation law present in the system \eqref{hughesD}  and then we will arrive to the well-posedness of  the discrete system.
%%%%%%
%    %
%%%%%%
\subsection{The eikonal equation on graph}\label{secHJ}
We study the graph eikonal equation  (see \cite{bm}, \cite{mos} for related results)
\begin{equation}\label{HJ}
\left\{
\begin{array}{ll}
    \max\limits_{y\sim x}\left\{-\frac{u^n(y)-u^n(x)}{w_{yx}}- \frac{1}{1-\r^n(y)}\right\}=0,\quad&x\in V_0, \, y\in V, \, n \in \N,\\[6pt]
    u^n(x)=0,&x\in V_b, \, n\in\N.
\end{array}
\right.
\end{equation}
 We assume that, for any $n\in\N$,
\begin{equation}\label{hj1}
   0\le  \r^n(x)\le 1-\d, \qquad \forall x\in V_0,
\end{equation}
for some $\d>0$, hence there exists   a  constant $M$ such that
\begin{equation}\label{hj0}
1\le \frac{1}{1-\r^n(x)}\le M, \qquad \forall x\in V.
\end{equation}
%%%
\begin{Theorem}
Let us assume that the condition \eqref{hj1} holds. Then, for any $n\in\N$, there exists a  unique solution of \eqref{HJ},   given by the formula
 \begin{equation}\label{hj2}
  u^n(x)=\min\left\{ \sum_{k=0}^{N-1}  \frac{w_{x_{k+1}x_k}}{1-\r^n(x_{k+1})}:\, y\in V_b,\g\in\cG_{xy}\right\}
 \end{equation}
$($with the convention $\sum_{k=0}^{-1}=0)$.
\end{Theorem}
%%%
\begin{Proof}\par
\emph{Existence. }
The function $u^n$ in \eqref{hj2} is well defined (i.e., the minimum is achieved) since, for each fixed $x$, the set of admissible paths connecting $x$ to the boundary  $V_b$  is finite.
 We first show that for $x\in V_0$
\begin{equation}\label{hj3}
 \max_{y\sim x}\left\{-\frac{u^n(y)-u^n(x)}{w_{yx}}- \frac{1}{1-\r^n(y)} \right\}\le 0.
\end{equation}
Given $y\sim x$,   let    $z\in V_b$ and  $\g=\{x_0=y,x_1,\dots, x_N=z\}\in \cG_{yz}$ be  such that
 \[u^n(y)= \sum_{k=0}^{N-1}  \frac{w_{x_{k+1}x_k}}{1-\r^n(x_{k+1})}.\]
 Then, $\{x, x_0,\dots, x_N=z\}\in \cG_{xz}$ is a   path connecting  $x$ to $z\in V_b$ and therefore
 \begin{align*}
    -(u^n(y)-u^n(x))&\le- \sum_{k=0}^{N-1}  \frac{w_{x_{k+1}x_k}}{1-\r^n(x_{k+1})}+  \frac{w_{x_0x}}{1-\r^n(x_0)}+ \sum_{k=0}^{N-1}  \frac{w_{x_{k+1}x_k}}{1-\r^n(x_{k+1})} \\
   &= \frac{w_{x_0x}}{1-\r^n(x_0)}=\frac{w_{yx}}{1-\r^n(y)},
 \end{align*}
from which we can conclude that \eqref{hj3} holds.\\ %\par
 Let us show now that for $x\in V_0$
\begin{equation}\label{hj4}
 \max_{y\sim x}\left\{-\frac{u^n(y)-u^n(x)}{w_{yx}}- \frac{1}{1-\r^n(y)} \right\}\ge 0.
\end{equation}
Let  $z\in V_b$ and  $\g=\{x_0=x,x_1,\dots, x_N=z\}\in \cG_{x,z}$ be   such that
\[u^n(x)= \sum_{k=0}^{N-1}  \frac{w_{x_{k+1}x_k}}{1-\r^n(x_{k+1})}.\]
Since   $x_1\sim z$ and  $\{x_1,\dots, x_N\}\in\cG_{x_1z}$, we get
 \begin{align*}
  -(u^n(x_1)-u^n(x))& \ge -\sum_{k=1}^{N-1}  \frac{w_{x_{k+1}x_k}}{1-\r^n(x_{k+1})}+\sum_{k=0}^{N-1}  \frac{w_{x_{k+1}x_k}}{1-\r^n(x_{k+1})}=\frac{w_{x_{1}x}}{1-\r^n(x_{1})}
 \end{align*}
 and, therefore,
\[
\max_{y\sim x}\left\{-\frac{u^n(y)-u^n(x)}{w_{yx}}- \frac{1}{1-\r^n(y)} \right\}\ge  -(u^n(x_1)-u^n(x))-\frac{w_{x_{1}x}}{1-\r^n(x_{1})}\ge 0.
\]
Combining \eqref{hj3} and \eqref{hj4}, we get  \eqref{HJ}.\par
Note that the positivity of the cost $1/(1-\r^n)$ implies $u^n(x)> 0$ for any $x\in V_0$.
If  $x\in V_b$, considering the stationary  path $\g=\{x_0=x\} \in\cG_{xx}$  which gives a null cost, we have
\[0\le u^n(x)=\min\left\{ \sum_{k=0}^{N-1}  \frac{w_{x_{k+1}x_k}}{1-\r^n(x_{k+1})}:\, y\in V_b,\g\in\cG_{xy}\right\}
\le 0\]
and therefore $u^n(x)=0$.\par
%%%%%%%%%%%%%%%%%%%%%%%%%%%%%%%
\emph{Uniqueness. } Let $u^n$, $v^n$ be two solutions of \eqref{HJ} and in addition we assume that
  $\max\limits_{V}\{u^n-v^n\}=\d$, for a strictly positive $\d$. We define also $W :=\arg\max\limits_{V}\{u^n-v^n\}$ and $m :=\min\{v^n(x):\, x\in W\}$. \par
Let  $x\in W$ be  such that $v^n(x)=m$. Since $u^n(x)=v^n(x)=0$ for $x\in V_b$, then $x$ belongs to $V_0$. Let $z\sim x$ be such that
\[
\max_{y\sim x}\left\{-\frac{v^n(y)-v^n(x)}{w_{yx}}- \frac{1}{1-\r^n(y)}\right\}=  -\frac{v^n(z)-v^n(x)}{w_{zx}}- \frac{1}{1-\r^n(z)} \, .
\]
Hence,
\begin{eqnarray*}
-\frac{v^n(z)-v^n(x)}{w_{zx}}- \frac{1}{1-\r^n(z)}=0\ge -\frac{u^n(z)-u^n(x)}{w_{zx}}- \frac{1}{1-\r^n(z)}
\end{eqnarray*}
%and therefore
from which $u^n(z)-v^n(z)\ge u^n(x)-v^n(x)=\d$. It follows that  $z\in W$ and, by $ -w^{-1}_{zx}(v^n(z)-v^n(x))\ge \frac{1}{1-\r^n(z)}> 0$, we get $m=v^n(x)>v^n(z)$, which
% gives a contradiction to the definition of $m$.
is in contradiction with the definition of $m$.
\end{Proof}
%%%%%%%%
In the next proposition, we give some regularity properties of $u^n$.
\begin{Proposition}
Let $u^n$ be the solution of \eqref{HJ}. Then
\begin{align}
    &  d(x,y)\le  u^n (x) \le  M d(x,y)  &&\quad \forall x\in V, y\in V_b, \label{regD1}\\[6pt]
    & |u^n(y)-u^n(x)|\le M d(x,y)&&\quad \forall x,y\in V, x\sim y, \label{regD2}
\end{align}
where  $d$ and $M$  are defined in  \eqref{geo} and \eqref{hj0}, respectively.
 \end{Proposition}
\begin{Proof}
Let $x\in V$ be. Then, for any  $y\in V_b$  and for any  $\g=\{x_0=x,x_1,\dots, x_N=y\}\in \cG_{xy}$, by the  inequalities in \eqref{hj0} we have
\[ \sum_{k=0}^{N-1}  w_{x_{k+1}x_k}\le  \sum_{k=0}^{N-1}  \frac{w_{x_{k+1}x_k}}{1-\r^n(x_{k+1})}\le  M \sum_{k=0}^{N-1}  w_{x_{k+1}x_k}.\]
Therefore the bounds \eqref{regD1} follow immediately.\par
Let $x\sim y$ and $\g=\{x_0=y,x_1,\dots, x_N=z\}\in \cG_{yz}$ be, where $z\in V_b$ is an optimal path for $u^n(y)$.
Then, $\{x, x_0,\dots, x_N=z\}\in \cG_{xz}$ is a   path connecting  $x$ to $z\in V_b$ and, therefore,
 \begin{align*}
    u^n(x)-u^n(y)&\le    \frac{w_{yx}}{1-\r^n(y)}+ \sum_{k=0}^{N-1}  \frac{w_{x_{k+1}x_k}}{1-\r^n(x_{k+1})}  - \sum_{k=0}^{N-1}  \frac{w_{x_{k+1}x_k}}{1-\r^n(x_{k+1})} \\
&= \frac{w_{yx}}{1-\r^n(y)}\le Md(x,y)
 \end{align*}
which proves the property \eqref{regD2}.
\end{Proof}
% Consider the problem
%\begin{equation}\label{eikD}
%\left\{
%\begin{array}{ll}
%    \max_{y\sim x}\left\{-\frac{u(y)-u(x)}{w_{yx}}- c(y)\right\}=0\quad&x\in V_0\\[4pt]
%    u(x)=0&x\in V_b
%\end{array}
%\right.
%\end{equation}
%where $c$ is a positive function on $V$.
%Define the function
Let us define the following function on the graph $\Gamma$:
\[
d_n (x,y):=\min\left\{ \sum_{k=0}^{N-1}  \frac{w_{x_{k+1}x_k}}{1-\r^n(x_{k+1})}:\, \g\in\cG_{xy}\right\}, \quad x,y\in V, n\in\N.
\]
Then, the solution of \eqref{HJ} can be written as
\begin{equation}\label{dist_boundaryD}
 u^n(x)=\inf\{d_n(x,y):\, y\in V_b\}
 \end{equation}
 (cf. with the formula \eqref{dist_boundary} in the continuous case).
Therefore,  $u^n$  is the distance from the boundary  taking into
account the distribution of the population on the graph:   the term  $1/(1-\r^n(y))$, which is the cost of passing from $x$ to $y$,
penalizes the vertices adjacent to $x$ with high population density.  The term $\textrm{sgn}(u^n(y)-u^n(x))$ in \eqref{hughesD}, which can be seen as the normalized discrete gradient of $u^n$,
gives the direction of the minimizing path to the boundary.  Moreover, if $\r^n\equiv 0$, then $d_n$ coincides with the path distance $d$ defined in \eqref{geo}.
%%%%%%%%%%%
%         %
%%%%%%%%%%%
\subsection{The conservation law on   graph }\label{secCL}
In this section we study the problem
%\begin{align}
\begin{eqnarray}
\left\{
\begin{array}{ll}
 \r^{n+1}(x)= \r^{n}(x)-\sum\limits_{y\sim x}\l h^n_{yx}\d^n_{yx},&\quad   x\in V, n\in\N, \label{CLD}\\
 \r^0(x)=\bar\r(x) &\quad x\in V,n=0,			% \label{initial}
\end{array}
\right .
\end{eqnarray}
%  \end{align}
where $h^n_{yx}$ satisfies \eqref{flux0}-\eqref{flux1}, $\l$ is a positive constant and $\d^n_{yx}$ is
equal to 1  if the flux is directed from $y$ to $x$ and to $-1$ viceversa (for \eqref{hughesD},
$\d^n_{yx}=\textrm{sgn}(u^n(y)-u^n(x))$).\\
We rewrite equation  \eqref{CLD} as
\begin{equation}\label{flux2}
    \r^{n+1}(x)=G(\r^n(x),\{\r^n(y)\}_{y\in I(x)})
\end{equation}
for a map $G:\R\times\R^{|V|}\to\R$.
%%%
\begin{Proposition}\label{prop:flux1}
Let us assume
\begin{equation}
   D \l \|m\|_{L^\infty(0,1)} \le 1, \label{CFL}
\end{equation}
where $D$ is defined in  \eqref{maxdeg} and $m$ in \eqref{flux1}.
Then, the  map $G$ is monotone in $[0,1]$, i.e. if $\r^n(x), \z^n(x) \in [0,1]$ for all $x\in V$, $n\in\N$, then
\[\r^n(x)\le \z^n(x)\quad \forall x\in V\quad\Rightarrow \quad \r^{n+1}(x)\le \z^{n+1}(x)\quad \forall x\in V.\]
\end{Proposition}
\begin{Proof}
Observe that
\[
G(\r^n(x),\{\r^n(y)\}_{y\in I(x)})=\r^n(x)-\sum_{y\sim x \atop \d^n_{yx}=1}\l h(\r^n(y),\r^n(x))+\sum_{y\sim x \atop \d^n_{yx}=-1}\l h(\r^n(x),\r^n(y)).
\]
We first prove that $\partial G/\partial \r^n(y)\ge 0$ for $y\in I(x)$. This follows immediately  by   \eqref{flux1} and by  the identity
\begin{equation}\label{flux2a}
\frac{\partial G}{\partial \r^n(y)}=
\left\{
                             \begin{array}{ll}
                            - \l  h_v (\r^n(y), \r^n(x)),  & \hbox{if $\d^n_{yx}=1$,} \\[4pt]
                              \l h_u(\r^n(x),\r^n(y)),  & \hbox{if $\d^n_{yx}=-1$.}
                             \end{array}
                           \right.
\end{equation}
Moreover, by \eqref{CFL} we have
\begin{multline}\label{flux2b}
     \frac{\partial G}{\partial \r^n(x)}
           =1-\sum_{y\sim x \atop \d^n_{yx}=1}\l h_u(\r^n(y),\r^n(x))+\sum_{y\sim x \atop \d^n_{yx}=-1}\l h_ v(\r^n(x),\r^n(y))\\
             \ge1- \sum_{y\sim x \atop\d^n_{yx}=1}\l m_+(\r^n(x))+\sum_{y\sim x\atop\d^n_{yx}=-1}\l m_-(\r^n(x))\ge  1-D\l\|m\|_{L^\infty(0,1)} \ge 0
\end{multline}
and therefore $G$ is    increasing  in $\r^n(x)$.
\end{Proof}
%%%%%%
\begin{Proposition}\label{prop:flux2}
Let us assume that \eqref{CFL} holds and that $0\le \bar\r(x)\le 1$ $\forall x\in V$. Then,
\begin{itemize}
  \item[(i)]   $0\le \r^n(x)\le 1$ $\quad\forall x\in V$, $n\in \N$.
  \item[(ii)]  If  $0\le \bar\r(x)<1$ and  the inequality in \eqref{CFL}
              is strict, then $0\le \r^n(x)<1$ $\forall x\in V$, $n\in \N$.
  \item[(iii)] $\sum\limits_{x\in V} \r^n(x)= \sum\limits_{x\in V} \bar\r(x)$, $\qquad\forall n\in \N$.
  \item[(iv)] $\sum\limits_{x\in V} |\r^{n+1}(x)-\r^n(x)|\le \sum\limits_{x\in V}| \r^1(x)-\bar\r(x)|$, $\qquad\forall n\in \N$.
\end{itemize}
\end{Proposition}
\begin{Proof}
By using \eqref{flux0}, the monotonicity of the map $G$ and the following
\[
0= G(0,\{0\}_{y\in I(x)})\le G(\r^n(x),\{\r^n(y)\}_{y\in I(x)})\le G(1,\{1\}_{y\in I(x)})=1,
\]
it follows that
\[0\le \r^n(x)\le1, \quad  \forall n\in\N,\; \forall x\in V. \]
Hence,  {\it (i)} holds.\\
If $D\l\|m\|_{L^\infty(0,1)}<1$,  by \eqref{flux2b} the map $G$ is strictly increasing in $\r^n(x)$.
Moreover,  by \eqref{flux2a}, $G(1,\{\bar\r(y)\}_{y\in I(x)})$ is the sum of terms non decreasing in $\bar\r(y)$.
Hence, if  $0\le \bar\r(x)<1$, then
\begin{equation*}\label{hj2c}
 \r^1(x)=G(\bar\r(x),\{\bar\r(y)\}_{y\in I(x)})<G(1,\{\bar\r(y)\}_{y\in I(x)})\le G(1,\{1\}_{y\in I(x)})=1, \quad\forall x\in V
\end{equation*}
and iterating on $n\in\N$  we get {\it (ii)}.\\
To prove the equality in {\it (iii)} , we  observe that
\begin{equation}\label{flux3}
 h^n_{yx}\d^n_{yx}+h^n_{xy}\d^n_{xy}=0\qquad \forall x,y\in V, \,x\sim y.
\end{equation}
In fact, if $\d^n_{yx}=1$, then $\d^n_{xy}=-1$ and by \eqref{flux} we have
\[h^n_{yx}\d^n_{yx}+h^n_{xy}\d^n_{xy}=h(\r^n(y),\r^n(x))-h(\r^n(y),\r^n(x))=0.\]
We proceed similarly if $\d^n_{yx}=-1$.
Since for each $x\in V$ there is a corresponding node $y\in V$ for which \eqref{flux3} holds,
we immediately get
%\[
\begin{equation}\label{flux7}
   \sum_{x\in V}\r^{n+1}(x)=\sum_{x\in V}\r^{n}(x)-\sum_{x\in V }\sum_{y\sim x}h^n_{yx}\d^n_{yx}
    =\sum_{x\in V}\r^{n}(x).
\end{equation}
%\\
%&+\sum_{x\in V}\sum_{y\sim x\atop y\in V}h^n_{yx}\d^n_{yx}+
%    \sum_{x\in V_0}\sum_{y\sim x\atop y\in V_b}h^n_{yx} =\sum_{x\in V}\r^{n}(x)- \sum_{x\in V_0}\sum_{y\sim x\atop y\in V_b}h^n_{yx}\\
% &  \le \sum_{x\in V_0}\r^{n}(x)
%\]
%i.e.
%\begin{equation}\label{flux7}
%    \sum_{x\in V} G( \r^n(x),\{ \r^n \}_{y\in I(x)}) = \sum_{x\in V}\r^{n}(x).
%\end{equation}
Iterating the previous argument on $n\in\N$ we get {\it (iii)}.\\
To prove {\it (iv)}, we consider the case $n=1$ and we observe that
\begin{multline}
    \sum_{x\in V}|\r^{2}(x)-\r^1(x)|=\sum_{x\in V}(\r^{2}(x)-\r^1(x))^+\sum_{x\in V}(\r^{1}(x)-\r^2(x))^+\\
   = \sum_{x\in V}\Big(G(\r^{1}(x),\{\r^1 \}_{y\in I(x)})-G(\bar\r(x),\{\bar\r \}_{y\in I(x)})\Big)^+ \\
    +\sum_{x\in V}\Big(G(\bar\r(x),\{\bar\r \}_{y\in I(x)})
    -G(\r^{1}(x),\{\r^1 \}_{y\in I(x)})\Big)^+.
 \label{flux8}
\end{multline}
Moreover, by the monotonicity of $G$, see {\it (i)}, and the mass conservation in  \eqref{flux7}, we can write
\begin{align*}
   & \sum_{x\in V}\Big(G(\r^{1}(x),\{\r^1 \}_{y\in I(x)})-G(\bar\r(x),\{\bar\r \}_{y\in I(x)})\Big)^+\\
&\le    \sum_{x\in V} \Big(G(\r^{1}\vee \bar \r (x),\{\r^1\vee \bar \r \}_{y\in I(x)})-G( \bar\r(x),\{\bar\r \}_{y\in I(x)})\Big)^+ \\
&=  \sum_{x\in V} G(\r^{1}\vee \bar \r (x),\{\r^1\vee \bar \r \}_{y\in I(x)})-G( \bar\r(x),\{\bar\r \}_{y\in I(x)}) \\
&=    \sum_{x\in V}  (\r^{1}\vee \bar \r) (x)-  \bar\r(x) =\sum_{x\in V}  ( \r^{1}  (x)-  \bar\r(x)  )^+,
\end{align*}
and similarly
\[
\sum_{x\in V}\Big(G(\bar\r(x),\{\bar\r \}_{y\in I(x)})
    -G(\r^{1}(x),\{\r^1 \}_{y\in I(x)})\Big)^+\le \sum_{x\in V} \Big(  \bar \r  (x)-   \r^1(x) \Big)^+.
\]
By substituting the previous inequality in \eqref{flux8} we obtain
\begin{align*}
   \sum_{x\in V}|\r^{2}(x)-\r^1(x)|\le \sum_{x\in V}  ( \r^{1}  (x)-  \bar\r(x)  )^+ +  (  \bar \r  (x)-   \r^1(x)  )^+=
      \sum_{x\in V}|    \r^{1}  (x)-  \bar\r(x)   |
\end{align*}
and, iterating, we get  {\it (iv)}.
\end{Proof}
%Property $(iii)$ in Prop.\ref{prop:flux2} gives the  mass conservation   for \eqref{CLD}.
The term
\[ \sum_{x\in V_b}\r^{n}(x)\]
represents the cumulative distribution of the population which has already reached the exit at the time $n$.
Since if $x\in V_b$ then $\d_{yx}^n=0$ $\forall y\in V_b$ with $y\sim x$, there is no flow inside  the boundary.
The assumption \eqref{flux1},  which gives a bound on the maximal admissible  velocity of the flux,  is exploited in conjunction with the assumption \eqref{CFL} in order to get the monotonicity of the map $G$ introduced in \eqref{flux2}.
 This property  guarantees that pedestrians can move only of one vertex for unit time.  Hence, people on not adjacent vertices of the graph cannot  interact in a single time interval.
\begin{Remark}\label{dirichlet}
For the numerical simulation, we also consider a homogeneous Dirichlet boundary condition in place of the no-flux boundary
condition \eqref{no_flux_BC}. The corresponding  conservation law on the graph is
\begin{eqnarray*}
\left\{
\begin{array}{ll}
 \r^{n+1}(x)= \r^{n}(x)-\sum\limits_{y\sim x}\l h^n_{yx}\d^n_{yx},&\quad   x\in V_0, n\in\N, \\
 \r^n(x)=0, &\quad   x\in V_b, n\in\N, \\[3pt]
 \r^0(x)=\bar\r(x),  &\quad x\in V, n=0.			% \label{initial}
\end{array}
\right .
\end{eqnarray*}
If we denote with $\tilde \rho^n$ the solution of the previous problem and by $\rho^n$ the solution of \eqref{CLD}, by the monotonicity of the scheme $G$ it is immediate to see that $\tilde \rho^n\le  \rho^n$ for any $n\in\N$. Hence,  also $\tilde \rho^n$ satisfies properties (i) and (ii)   in Prop. \ref{prop:flux2}.
\end{Remark}
As an immediate consequence of the Proposition \ref{prop:flux2} and the assumption \eqref{hj1}, we have the well-posedness
of the Hughes' model on a graph.
\begin{Corollary}
Assume that $\l D\|m\|_{L^\infty(0,1)}<1$ and $0\le \bar\r(x)<1$ $\forall x\in V$. Then the problem \eqref{hughesD}
is well defined $\forall n\in\N$.
\end{Corollary}
\begin{Proof}
By Proposition \ref{prop:flux2}{\it (ii)} and the condition \eqref{hj1}, the eikonal equation \eqref{HJ} is well defined $\forall n\in \N$.
It follows that also the conservation law \eqref{CLD} is well defined $\forall n\in \N$.
\end{Proof}

%%%%%%%%%%%%
%          %
%%%%%%%%%%%%
\section{Numerical results for the  Hughes' model  on graph}\label{Sec4}
In this section we discuss the numerical implementation of   the discrete Hughes' system \eqref{hughesD}, which is  considered as a discretization
of the continuous Hughes' system \eqref{hughesnetwork} on a network $\cN$. We introduce a time step $dt$ and a spatial step $dx$.
The discretization of the edges of the network $\cN$    is done uniformly with respect to $dx$ with the points
of the discretization giving the vertices of the graph $\G$. The   stability  condition in \eqref{CFL} is verified if
$$  dt  \le \frac{  dx }{D \|m\|_{L^\infty(0,1)}}, $$
being $\lambda = dt/dx$ in \eqref{CLD}.
The choice of the flux function $h$ is an important point. It is well known that a  low order scheme gives poor accuracy on smooth regions of the solutions, conversely a  high order scheme could develop spurious oscillations    bringing to a  high local error in non smooth regions of the solutions.  We experimentally observed that a good compromise is represented by  the Engquist-Osher flux \eqref{EO}.
%\begin{equation}\label{ENO}
%      h(\r^n(y),\r^n(x))=\half\Big(g(\r^n(y))+g(\r^n(x)))\Big)+\half\int_{\r^n(x)}^{\r^n(y)} |g'(\r)|\ud\r.
%\end{equation}
For this reason, from now on we will consider this form of the numerical flux $h$.
%In the rest of the section we consider   this form of the numerical flux $h$ with the additional second order term (see pp. 7) where indicated.\\
In the resolution of  the system \eqref{hughesD},   the conservation law  is
explicit in time,  hence its computation does not present any difficulty. The resolution of  the eikonal equation
is more delicate and we consider a \emph{value iteration} technique.  Taken an initial guess $u_0(x)$, we iterate   the explicit system
\begin{eqnarray*}
\left\{
\begin{array}{ll}
 v^{k+1}(x)= \min\limits_{y\sim x}\left\{ v^k(y)+\frac{w_{yx}}{1-\r^n(y)}\right\}, & \quad x\in V, \\
 v^{k+1}(x)= 0, & \quad x\in V_b, \\
   v^{0}(x)=u_0(x). &  \
\end{array}
\right .
\end{eqnarray*}
Under  some non restrictive hypotheses (see \cite{PB79}), such iteration is a contraction and  converges monotonically for  $k\rightarrow+\infty$ to the solution $u^n$ of  \eqref{HJ}.  It has been   shown the relevance of a good initial guess $u_0(x)$ to have a fast convergence (cf. \cite{KAF13}). A perfect candidate to   play this role  is the function $u^{n-1}(x)$ that is the  solution of the discrete eikonal equation at the previous time step.
%
%
%
%We summarize the technique proposed in the following pseudo code
%
%
%
%%\begin{table}[h]
%\begin{center}
%\line(1,0){435}\\
%\vspace{0.2cm}
%{\sc Discrete Hughes Algorithm}
%\vspace{-0.5cm}
%\line(1,0){435}
%\end{center}
%\vspace{0.3cm}
%\phantom{e}
%Given an initial density $\bar\r(x)$, set $\r^0(x)=\bar\r(x)$ for every $ x\in V$. \par
%\phantom{aaaaa} Set the stopping parameter the time step $dt$ and $\varepsilon$ to a desired value.\\
%\phantom{aaaaa} Initialize $u_0(x)$ for $x\in V$. \\
%\phantom{aaaaa} Set $n:=0$;
%\begin{enumerate}
%\item[1)] \emph{(Computation discrete eikonal equation)}\par
%Set $k:= 0$,
%\begin{itemize}
%\item[1.i)]  Set  $v^{0}(x)=u_0(x)$
%\item[1.ii)] Compute $ v^{k+1}(x)= \min\limits_{y\sim x}\left\{ v^k(y)+\frac{w_{yx}}{1-\r^n(y)}\right\}$ for every $x\in V$.
%\item[1.iii)] If $\max_{x\in V}|v^{k+1}(x)-v^{k}(x)|\leq \varepsilon$ \\
% \vspace{0.1cm}
% \phantom{aaaaa}then set $u^n(x):=v^{k+1}(x)$ for every $x\in V$ and go to 2), \\
% \vspace{0.1cm}
%\phantom{aaaaa}otherwise set $k:=k+1$ and go to 1.i)
%\end{itemize}
%\item[2)]  \emph{(Time evolution of the conservation law)}
%\begin{itemize}
%\item[2i)] Compute $\r^{n+1}(x)= \r^{n}(x)-\sum\limits_{y\sim x}\frac{dt}{w_{yx}} h^n_{yx}\textrm{sgn}(u^n(y)-u^n(x))$ for each $x\in V$
%\item[2ii)] Set $t:=t+1$, $u_0(x)=u^n(x)$ for every $x\in V$ and go to 1)
%\end{itemize}
%\end{enumerate}
%\vspace{-0.2cm}
%\line(1,0){435}\\
%\vspace{0.4cm}
%%\caption{Pseudo-code of the technique to solve \eqref{hughesD}}\label{Al}
%%\end{table}

\subsection{Synthetic Tests}
\begin{figure}[ht]
\begin{center}
\includegraphics[height=4.5cm]{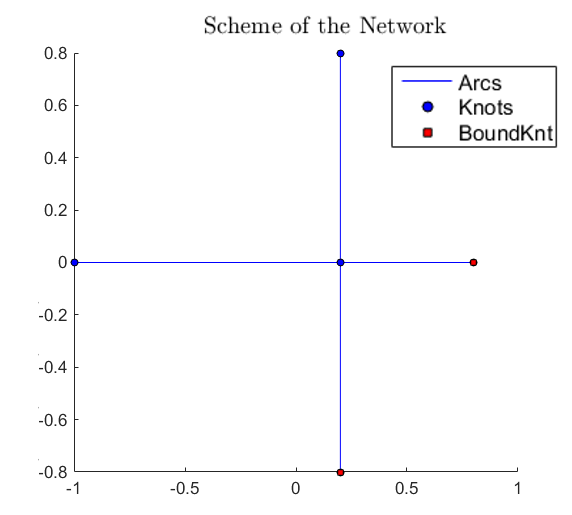} %[height=5.2cm]
\hspace{0.3cm}
\includegraphics[height=4.8cm]{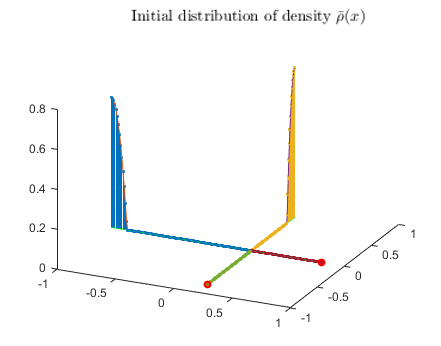}
\caption{Scheme of the network and initial density.} \label{fig1}
%\caption{Test 1: scheme of the graph and initial density.} \label{fig1}
\end{center}
\end{figure}

In this section we consider a simple network composed  of five nodes and four edges (see Figure \ref{fig1}, left).

\begin{figure}[ht]
\begin{center}
\includegraphics[height=7cm]{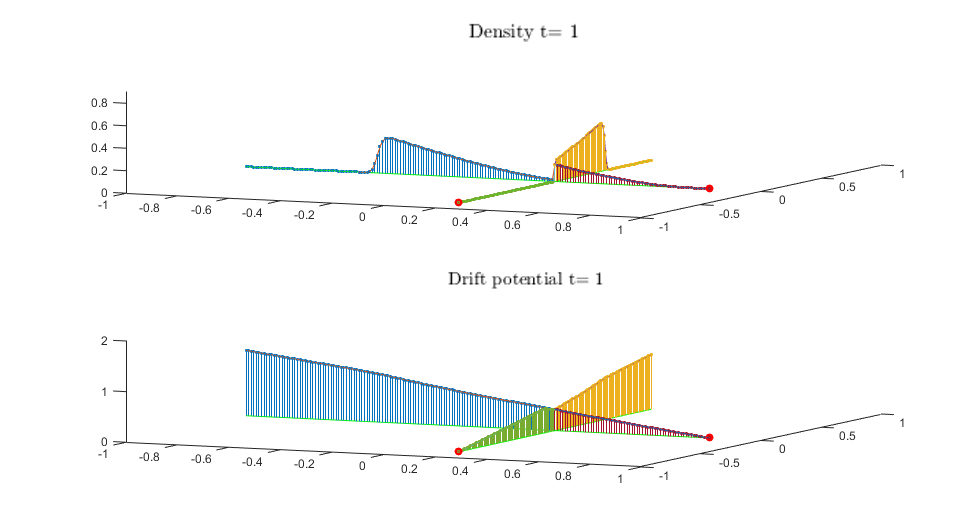} %[height=5.2cm]
\caption{Test 1 (Dirichlet boundary conditions): density and potential before the first time of  interaction.} \label{fig2}
%\caption{Test 1: density and potential before the first time of  interaction.} \label{fig2}
\end{center}
\end{figure}

The initial density $\bar \r(x)$ (see Figure \ref{fig1}, right) is defined by the restriction on the network of the function
$$  \bar \r(x_1,x_2) := \max(0, 0.65-  4\,(x_1+1)^2-4\,x_2^2,\, 0.75-  (6\,(x_1-0.2))^2-(6\,(x_2-0.8))^2).$$
The set of the boundary points  $V_b$, i.e. the \emph{target points}, is given by the vertexes in $(0.2,-0.8)$ and $(0.8,0)$, where   $u^n$ is imposed equal to zero.

We will consider two possible cases for the boundary conditions (BCs) for the conservation law: the case with a no-flux   condition  (in such case the mass is conserved inside the network) and the case with a homogeneous Dirichlet condition  on the target points  (see Remark \ref{dirichlet}).
Those BCs are relative to a different choice of the model: the no-flux condition corresponds to target points where the crowd tends to concentrate,
for example the stage of a concert, the various points of interest during the annual Hajj, see \cite{anon}, etc. The homogeneous Dirichlet condition, instead, corresponds to    target points which can be seen as \emph{exits} of large dimensions: any mass touching them exits instantaneously from the network.\\
First objective of this section is to show the stability of the discrete system: with this aim we consider a first order numerical flux as in \eqref{EO}. The case with a second order correction (stochastic perturbation adding diffusion) is more regular and it will be take into account in the next Test 3.

%From the modellistic point of view this request correspond to the imposition of some Dirichlet boundary conditions on the continuous system \eqref{Dirichlet_BC}, typical description of an , which is at the same time the objective of the agents of the system and a point where the mass can exit from the domain.

\noindent We perform the simulation fixing the discretization parameter $dx=0.01$ and the time step $dt=0.002$, hence the stability condition \eqref{CFL} is   verified since we observe $\|m\|_\infty \leq 1-2\rho\leq 1$ and $dx/(D\,dt)=5/4$.

\begin{figure}%[ht]
\begin{center}
\includegraphics[height=6cm]{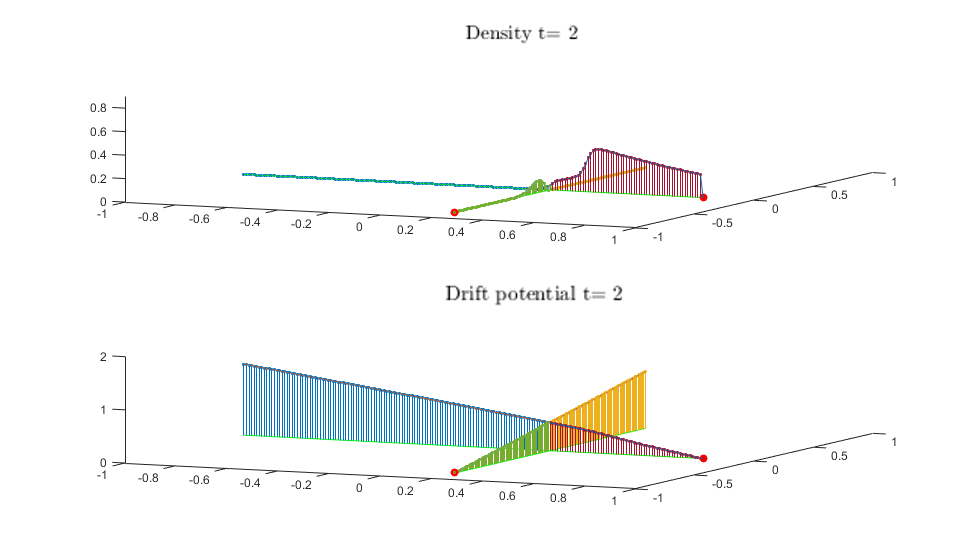} %[height=5.2cm]
\caption{Test 1 (Dirichlet boundary conditions): density and potential after  the first time of  interaction  at  $(0.2,-0.8)$.} \label{fig3}
%\vspace{1cm}
\end{center}
\end{figure}

%%%%%%% TEST 1 DIRICHLET BCs
\paragraph{\bf Test 1}
We start considering the case of homogeneous Dirichlet boundary conditions. At the beginning of the simulation, the two initial masses start to move in the direction of the two target points acting as exits.   The mass coming from the edge connecting $(0.2,0)$ to $(0.2,0.8)$ has a shorter path   so it arrives on the junction node and it turns in direction of the closer exit $(0.8,0)$ (Figure \ref{fig1}). In Figure \ref{fig2} it is possible to observe the first interaction time between the two masses:   the arrival of the   mass coming from the edge connecting $(-1,0)$ to $(0.2,0.8)$  produces a congestion near the exit $(0.8,0)$, therefore the other exit $(0.2,-0.8)$ becomes convenient  as it is possible to observe in the graph of the drift potential $u$, see Figure \ref{fig3} (below) and compare with Figure \ref{fig2} (below). Therefore,  the exit $(0.2,-0.8)$  attracts a part of the mass (Figure  \ref{fig3} above). This phenomena is peculiar of the Hughes' model and it has been observed also in other works (see e.g.  \cite{h03} or  \cite{dmpw}). Once reached the target, the mass exits from the network (Figure  \ref{fig3} above, on the exit $(0.8,0)$).

\begin{figure}%[ht]
\begin{center}
\includegraphics[height=3.5cm]{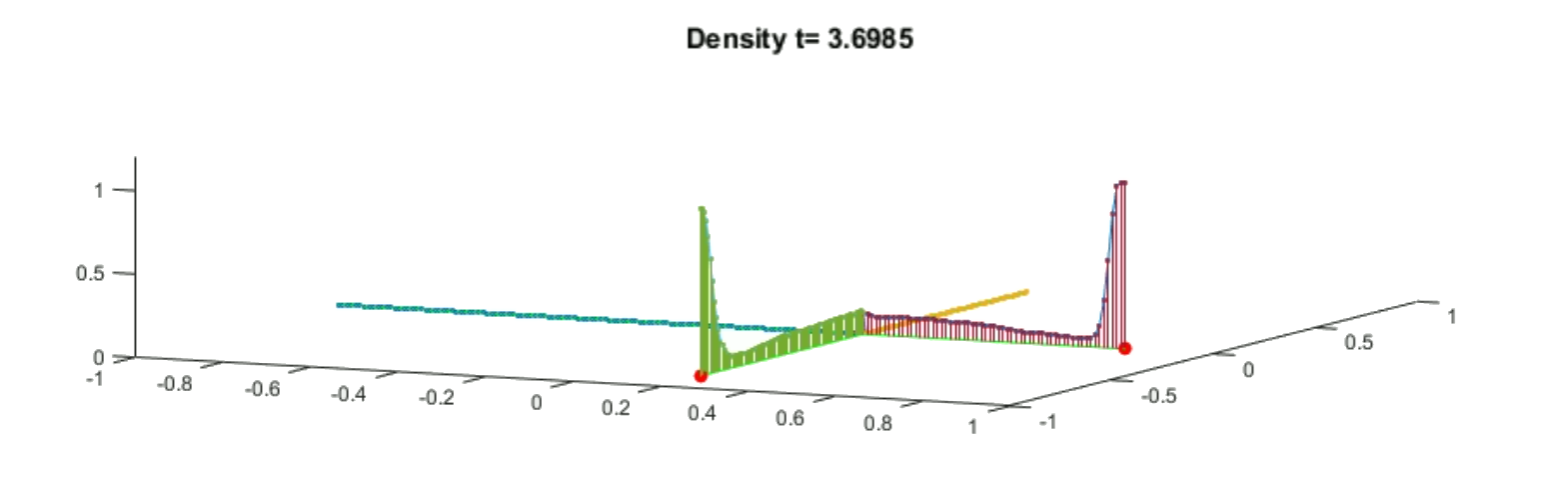} %[height=5.2cm]
\caption{Test 2 (No-flux boundary conditions): stable configuration obtained for $t>3.5$.} \label{fig32}
%\caption{Test 1 (No boundary conditions): stable configuration obtained for $t>3.5$.} \label{fig32}
\end{center}
\end{figure}

%%%%%%% TEST 2 NO BCs
\paragraph{\bf Test 2}
In a second simulation we compute the same solution with the no-flux  boundary condition. In this case, the mass is conserved. The first part of the test shows the same results as above: the masses are attracted by the target point $(0.8,0)$, but differently the mass starts to concentrate, reducing its speed until the exit $(0.2,-0.8)$ becomes convenient. When also this second target point reaches its maximal value of the density getting congested, the mass reaches a stable configuration  (Figure \ref{fig32}). In the case of a coarse time step $dt$, some oscillatory phenomena are observed in the second part of the test (chattering), where essentially the mass changes alternatively objective between the two target points. With a finer time discretization, those effects are reduced till disappearing. This test confirms even more than the previous case the stability of the discrete system and the fact that the density is always strictly lower  than the maximal value 1, even in the extreme case, when we force the mass to concentrate.%\\

\begin{figure}[ht]
\begin{center}
\includegraphics[height=6.8cm]{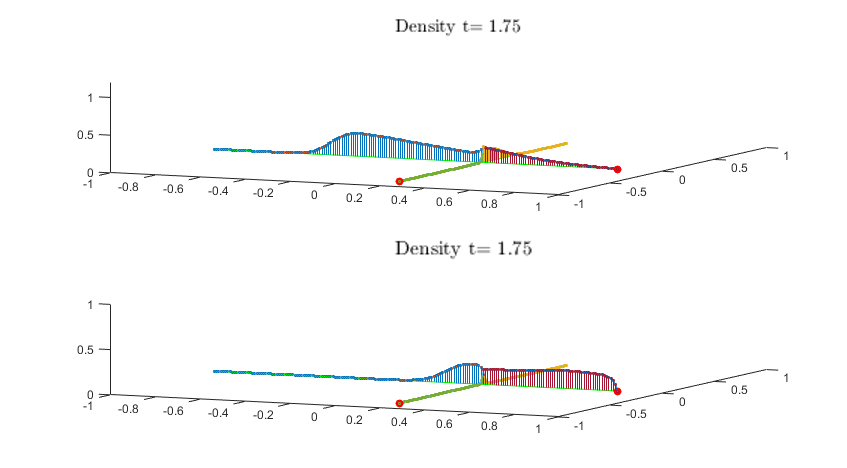} %[height=5.2cm]
\caption{Test 3 (Dirichlet BCs with diffusion $\epsilon=1$): Density at two different time steps ($t=0.75$ and $t=1.75$).} \label{fig33}
%, showing that congestion is not present.
\end{center}
\end{figure}

%%%%%%% TEST 3 Dirichlet BCs + diffusive term
\paragraph{\bf Test 3}
As last simulation on this network, we add to the conservation law   a term of the type \eqref{diff}, which can be interpreted, from a model point of view, as a stochastic perturbation in the motion of the mass and, from an analytic point of view, as a second order regularizing term in the equation. In Figure \ref{fig33} it is possible to see the effects of the diffusive term: the solution is more regular and congestion is not present. This is a phenomenon observed also in \cite{CFSW}: the presence of a stochastic noise prevents the mass to concentrate over a certain ratio. This has the indirect effect to help the overall \emph{evacuation time} (i.e. the first time step where the density on the domain is null everywhere) for certain configurations of the system (we can observe this comparing Figure \ref{fig3} with Figure \ref{fig33}). Avoiding congestion brings also some other macroscopic effects: in this case all the mass is exiting by the more convenient ``exit'' located in $(0.8,0)$. Since this arc is not getting congested, the agents do not change strategy using the other ``exit'' in $(0.2,-0.8)$.

\subsection{A Stadium evacuation test }
\begin{figure}[t]
\begin{center}
\includegraphics[height=5cm]{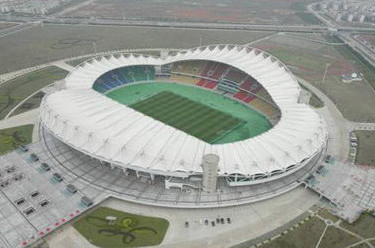}
\includegraphics[height=7cm]{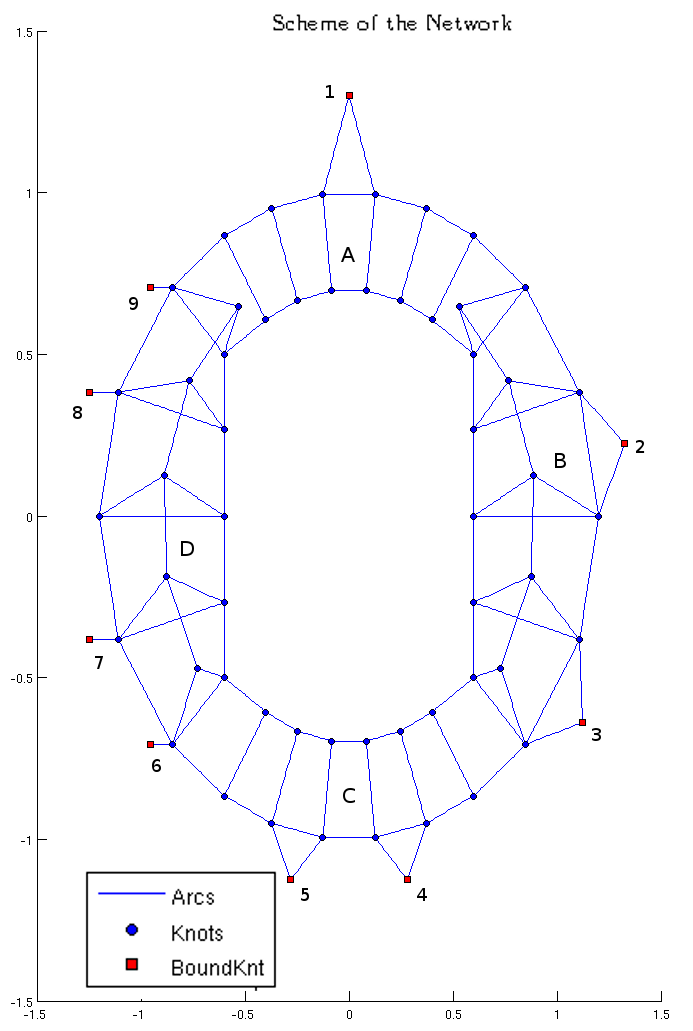}
\caption{The Wuhan Sports Centre (left) and the evacuation network considered in our study (right).} \label{stad}
\end{center}
\end{figure}

The stadium at the Wuhan Sports Centre (Fig. \ref{stad}, left) is a multi-use stadium located in Wuhan, China. Completed in 2002, it was used as test benchmark for mass-evacuation in \cite{Fang}. The
stadium has 42 bleachers (tiers of seats) distributed on all 3 floors
and has 9 exits (Fig. \ref{stad}, right) for evacuation; the capacity declared of the structure is of 54,357 spectators. Transforming a bit the structure (we consider a planar network) the evacuation network in
this stadium (Fig.  \ref{stad}, right) has 108 arcs and 63 nodes. After an uniform discretization of the arcs, the number of nodes of the network is around $7\cdot 10^3$ with a similar number of connections.

%In the second test we consider the more complex network represented in Figure \ref{fig4}. It is close to a simplified version of structures and devices in order to slow the speed of a density (crowd, information) and reduce the rate of exit from a system (non capisco bene cosa si intende). (add references from Engineering journals). The  initial density $\bar \r(x)$ is given by the restriction to the network of the function
%$$  \bar \r(x_1,x_2):= \max(0, 0.65-  (x_1+1)^2-x_2^2),$$
%see Figure \ref{fig4} (right),     $V_b$ is given by the vertex $(0.8,0)$ and the approximation parameters are set to $dx=0.01$ and $dt=0.005$.\\
%In Figure \ref{fig5}, we see the density $\r$ at the times   $t=1$ and $t=2$. Since the closest way to reach the exit is  through the central edge, this is the first one to be filled (Figure \ref{fig5} above). After congestion of the central edge, the choice of the other edges becomes progressively more convenient and they  are chosen by a part of the population (Figure \ref{fig5} below). The symmetry of the scheme can help to check the correctness of the approximation: indeed  the evolution of the mass is symmetric with respect to the central edge. \textbf{plotterei anche il potenziale a $t=1$ and $t=2$ }

\begin{figure}[t]
\begin{center}
\includegraphics[height=7cm]{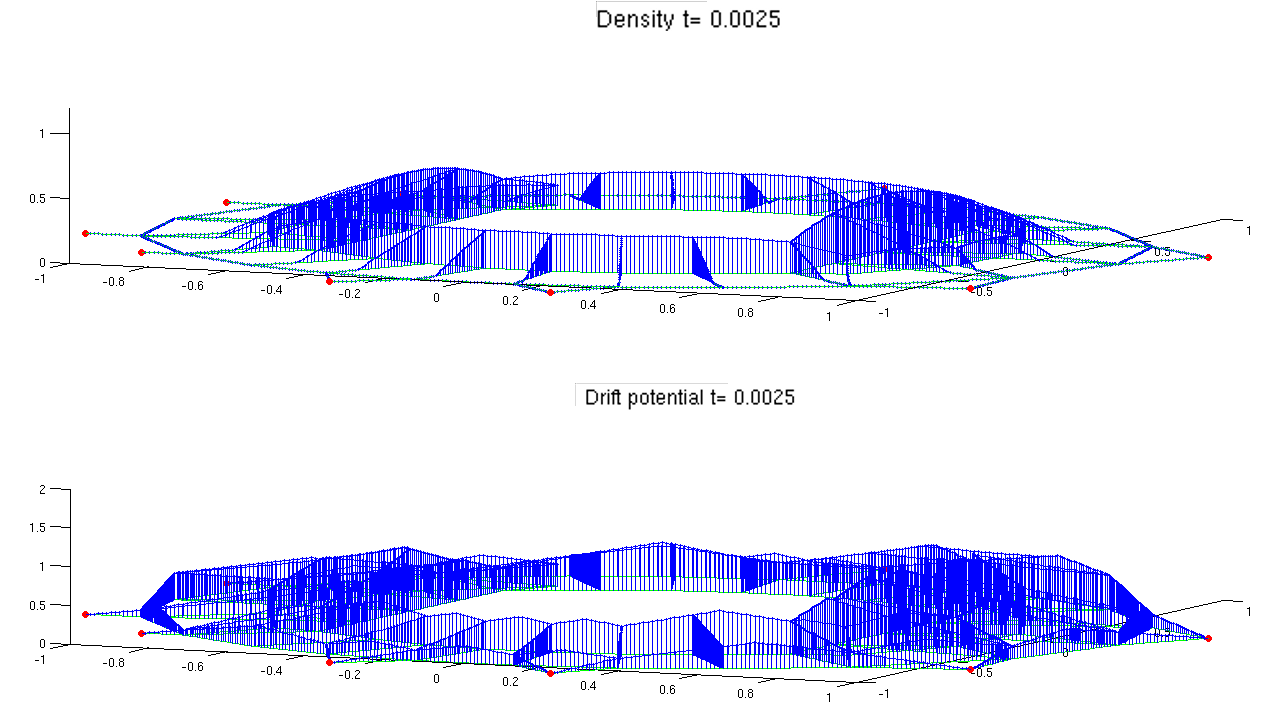}
\caption{Initial distribution of density on the network (up) and drift potential in the initial configuration (down).} \label{stad1}
\end{center}
\end{figure}

\begin{figure}[t]
\begin{center}
\includegraphics[height=7.2cm]{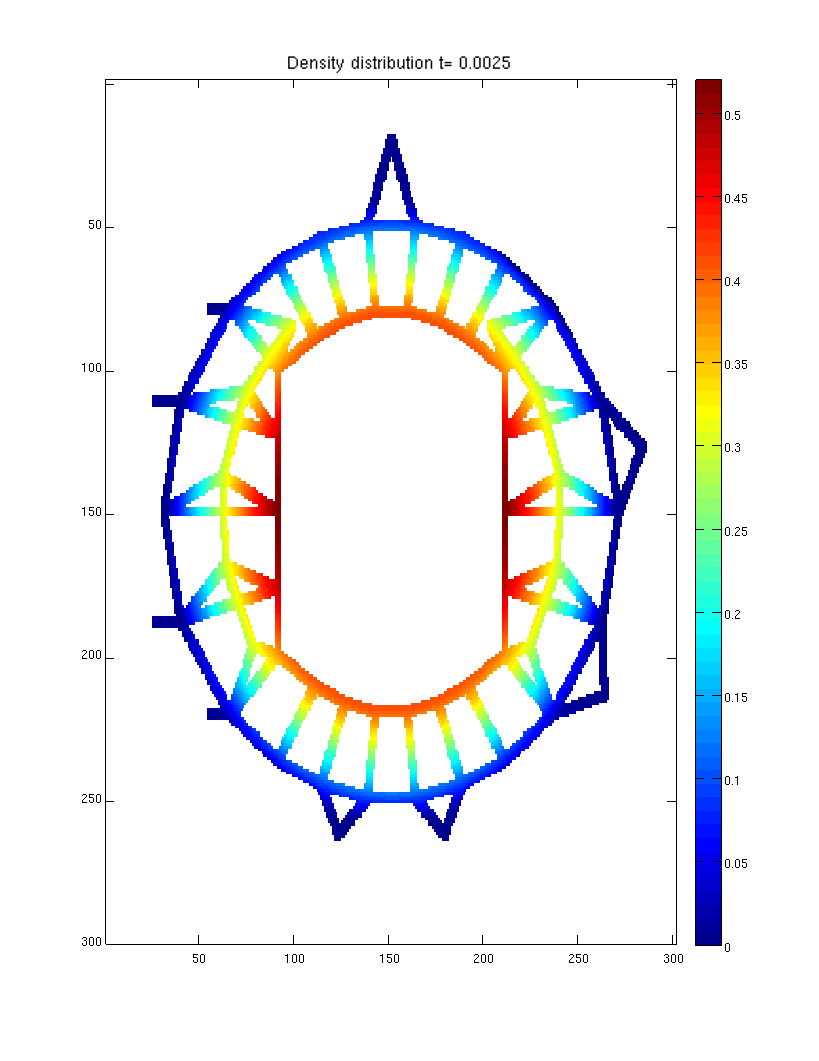}
\includegraphics[height=7.2cm]{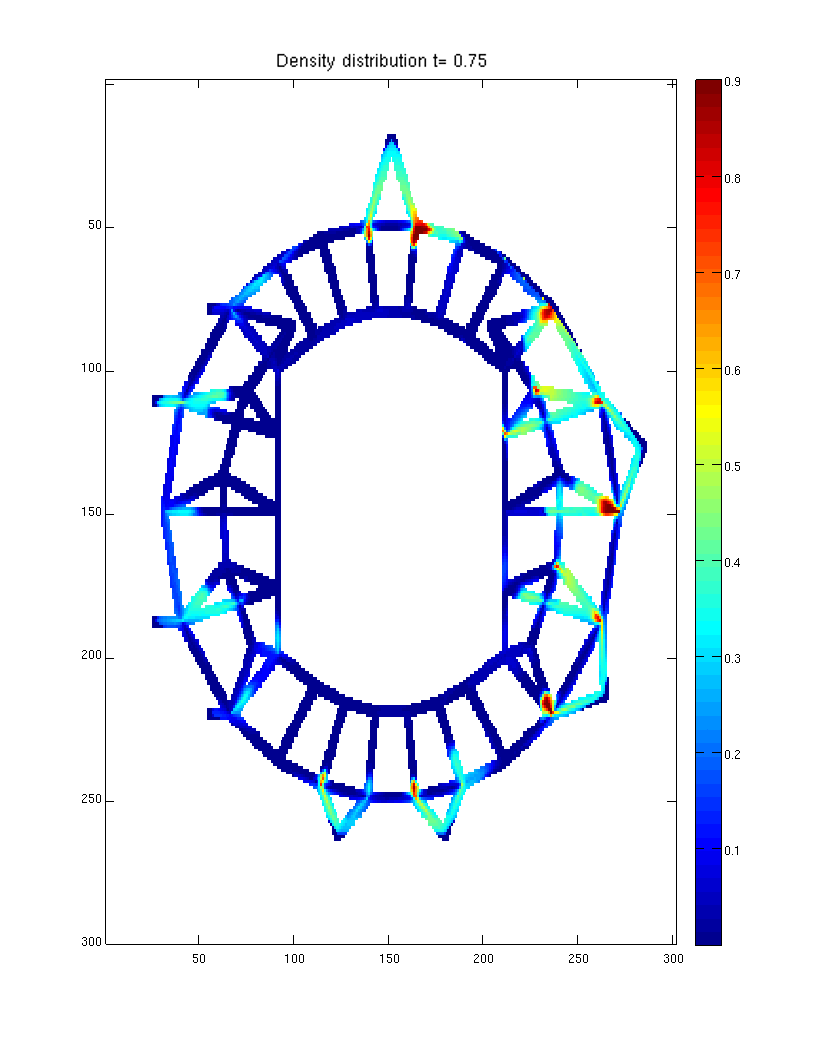}
\includegraphics[height=7.2cm]{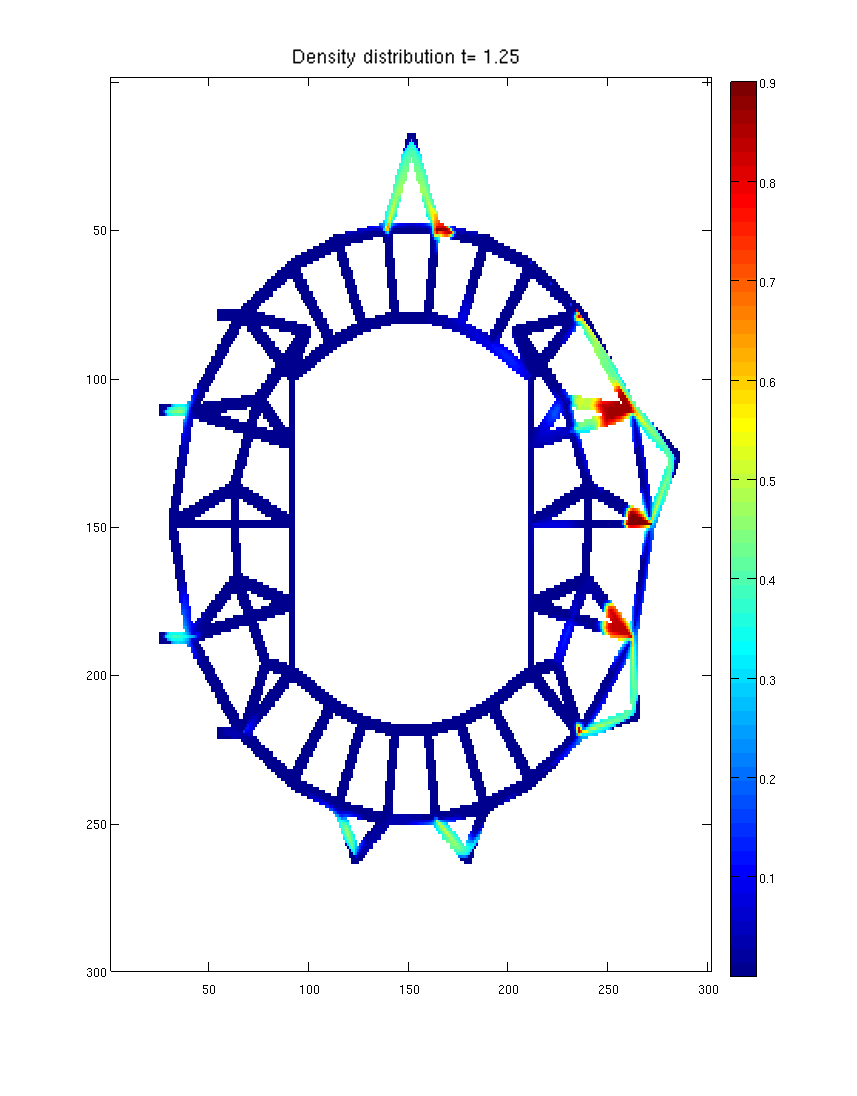}
\includegraphics[height=7.2cm]{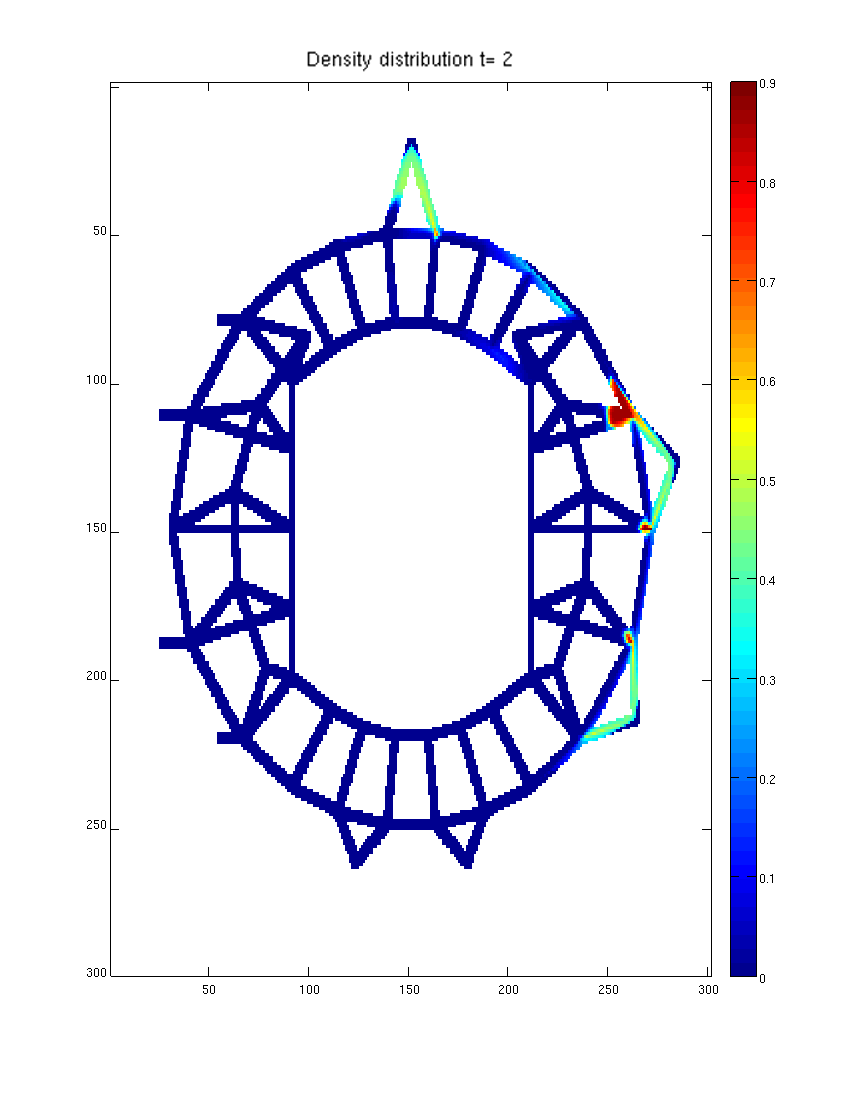}
\caption{Distribution of density on the network at various moments of the evolutions. Respectively (from left to right, up to down) $t=0.0025$, $0.75$, $1.25$, $2$.} \label{stad3}
\end{center}
\end{figure}

The choice of the initial configuration of density can be variable with respect to the aspect that we want to underline (by choosing a high initial uniform density distribution, we can test the network in an extremely crowded scenario; a random density choice can simulate some not standard cases of anomalous local concentration of crowd, etc). In this test, we chose the following initial datum
$$  \bar \r(x_1,x_2) := \max(0,0.7-0.7 x^2-0.84 y^2),$$
(we always mean the restriction of such function on the network), this
distribution in our intention should render the higher concentration of spectators in the areas closer to the court. We approximate uniformly the arcs using the discretization step $dx=0.01$ and we sample the time with $dt=0.002$. Also in this case the condition  \eqref{CFL} is guaranteed (the maximum number of connections per node $D$ is $5$) and the scheme is stable.
In Figure \ref{stad1} we can see the initial distribution of the density and the potential driving the individuals toward one of the exits. In Figure \ref{stad3} it is shown the evolution of the system in various moments.  We can notice that, despite the general symmetry of the structure, the behavior is highly conditioned by the position and the number of the exits. Of particular interest is the difference between the evacuation of sectors A/B and C/D (refer to Fig. \ref{stad}). As it can be seen from the scheme, the sector A is served by only one exit, differently from C, where two exits are present. Analogously, the sector D has 4 exits conversely to sector B which has just two. This brings to an orderly and efficient evacuation in the sectors C/D, with a well balanced use of the exits available. In the other case (sector A/B), the observed dynamics are different: on the paths toward the only three exits, in proximity to some nodes with multiple access, there are the appearance of high density congested regions.  We can observe also some of the phenomena discussed previously: reduction of the speed in the congested regions, changes of strategy, doubtful choice between two strategies (chattering). This has the macroscopic effect to rise up the final time necessary to evacuate the regions involved as we can observe in the last samples of Figure \ref{stad3}: the sectors C/D are already empty, the sectors A/B, instead, show congestion and a laborious flow through them. It is interesting to observe that a more efficient (for evacuation purposes) network is not trivially a network with more exits but a structure that avoids to drive big masses of agents to pass at the same time in the same nodes. This is a consequence to the incapacity of the agents represented in the model to forecast the future configuration of the system in order to choose the best strategy to adopt. For those reasons, the model is particularly appropriate to simulate the behavior of a crowd in a known network in presence of unpredicted events (an evacuation order, unusual high concentration in common transport facilities, etc.).

%%%%%%%%%%%%
%          CONCLUSIONS %
%%%%%%%%%%%%
\section{Conclusions}\label{Sec5}
In this paper we have presented a discrete Hughes'  model for pedestrian flow on a graph.
We have shown that, differently from the analogous continuous model, this discrete model is well-posed for any time $n\in \N$ under some natural assumptions on the flux, continuing  to share some qualitative properties  with the corresponding continuous model, as the interpretation of the solution of the graph eikonal equation as a distance from the boundary, change  of strategies, congestion, etc. \\
Several tests have been shown, analyzing and comparing the results and the behaviors obtained with different conditions (no BCs, homogeneous Dirichlet BCs or adding a diffusive term).
The experimental examples have confirmed the validity of the proposed model,  showing that the discrete system is always stable, even in the extreme case, when we force the mass to concentrate.
%The last test related to the evacuation of a stadium

%%%%%%%%%%%%
%          %
%%%%%%%%%%%%

\end{document}